\documentclass[secthm,seceqn,amsthm,ussrhead,10pt]{amsart}
\usepackage{amsmath,latexsym}
\usepackage[english]{babel}
\usepackage[psamsfonts]{amssymb}
\usepackage{times}
\usepackage{cite}
\usepackage{pdflscape} 
\usepackage{ulem}
\usepackage[mathcal]{euscript}
\usepackage{tikz}
\usepackage{hyperref}
\usepackage{cancel}
\usetikzlibrary{arrows}

\newtheorem{Th}{Theorem}
\newtheorem{Lem}[Th]{Lemma}

\newtheorem{corollary}[Th]{Corollary}

\newenvironment{Proof}[1][Proof.]{\begin{trivlist}
\item[\hskip \labelsep {\bfseries #1}]}{\flushright
$\Box$\end{trivlist}}

\begin{document}
	\sloppy


{\Large Degenerations of  Zinbiel and nilpotent Leibniz algebras
\footnote{The work was supported by RFBR 14-01-00014, 16-31-00004 and 16-31-50017; by R\& D 6.38.191.2014 of Saint-Petersburg State University, "Structure theory, classification, geometry, K-theory and arithmetics of algebraic groups and related structures"; and by the President's Program "Support of Young Russian Scientists" (grant MK-1378.2017.1).}}

\medskip

\medskip

\medskip

\medskip
\textbf{Ivan Kaygorodov$^{a}$, Yury Popov$^{b,c}$, Alexandre Pozhidaev$^{c,d},$ Yury Volkov$^{e}$}
\medskip

{\tiny
$^{a}$ Universidade Federal do ABC, CMCC, Santo Andr\'{e}, Brazil.

$^{b}$ Universidade Estadual de Campinas, Campinas, Brazil.

$^{c}$ Novosibirsk State University, Novosibirsk, Russia.

$^{d}$ Sobolev Institute of Mathematics, Novosibirsk, Russia.

$^{e}$ Saint Petersburg State University, Saint Petersburg, Russia.
\smallskip

    E-mail addresses:\smallskip

    Ivan Kaygorodov (kaygorodov.ivan@gmail.com),
    
    Yury Popov (yuri.ppv@gmail.com),

    Alexandre Pozhidaev (app@math.nsc.ru),
    
    Yury Volkov (wolf86\_666@list.ru).

}

       \vspace{0.3cm}

{\bf Abstract.} 
We describe degenerations of four-dimensional Zinbiel and four-dimensional nilpotent Leibniz algebras over $\mathbb{C}.$  
In particular, we describe all irreducible components in the corresponding varieties.\smallskip

       \vspace{0.3cm}

       \vspace{0.3cm}

\section{Introduction}

       \vspace{0.3cm}

Degenerations of algebras is an interesting subject, which was studied in various papers (see, for example, \cite{B99,B05,M79,M80,AOR05,CKLO13,BC99,S90,GRH,GRH2,BB09,BB14}).
In particular, there are many results concerning degenerations of algebras of low dimensions in a  variety defined by a set of identities.
One of important problems in this direction is the description of so-called rigid algebras. These algebras are of big interest, since the closures of their orbits under the action of generalized linear group form irreducible components of a variety under consideration
(with respect to the Zariski topology). For example, the rigid algebras were classified in the varieties of
low dimensional associative (see \cite{M79,M80}) and Leibniz (see \cite{AOR05}) algebras.
There are fewer works in which the full information about degenerations was found for some variety of algebras.
This problem was solved for four-dimensional Lie algebras in \cite{BC99}, for nilpotent five- and six-dimensional Lie algebras in \cite{S90,GRH}, 
for three-dimensional nilpotent Leibniz algebras in \cite{rom}, 
for  two-dimensional pre-Lie algebras in \cite{BB09},  for three-dimensional Novikov algebras in \cite{BB14},  
for nilpotent five- and six-dimensional Malcev algebras in \cite{kpv}, 
and all $2$-dimensional algebras \cite{kv16}.

The Leibniz algebras were introduced by Bloh in \cite{bloh} as a generalization of Lie algebras.
The study of the structure theory and other properties of Leibniz algebras was initiated by Loday in \cite{lodaypir}.
Leibniz algebras were also studied in \cite{leib1,leib2,yau,leib4}.
Zinbiel algebras were introduced by Loday in \cite{loday}.
Under the Koszul duality the operad of Zinbiel algebras is dual to the operad of Leibniz algebras.
 Zinbiel algebras were studied in \cite{abash, AKO10,cam13, dok, dzhuma5}.

In this paper we give the full information about degenerations of Zinbiel and nilpotent Leibniz algebras of dimension $4.$
More precisely, we construct the graph of primary degenerations. 
The vertices of this graph are the isomorphism classes of algebras in the variety under consideration.
An algebra $A$ degenerates to an algebra $B$ iff there is a path from the vertex corresponding to $A$ to the vertex corresponding to $B$.
Also we describe rigid algebras and irreducible components in the varieties of algebras under consideration.

{\bf Remark.} 
We give the classification of rigid algebras and irreducible components in the variety of nilpotent
four-dimensional Leibniz algebras in Corollary 6. One can note that this classification slightly differs
from the one obtained in \cite{AOR05}. It happens due to the following reasons:
\begin{enumerate}
\item The algebras $\mathfrak{L}_i \ (9 \leq i \leq 12)$ are missed in the classification of nilpotent four-dimensional Leibniz
algebras used in \cite{AOR05}. Thus, the rigid algebra $\mathfrak{L}_{11}$  is replaced by the nonrigid algebra $\mathfrak{L}_1$ there. Note
also that $\mathfrak{L}_{11}$ degenerates to $\mathfrak{N}_5$ and kills its rigidity in the variety of nilpotent Leibniz algebras, but,
as it will be shown in this paper, the algebra $\mathfrak{N}_5$ is nonrigid even in the variety of three-step nilpotent
algebras.
\item In the classification of nilpotent four-dimensional Leibniz algebras used in \cite{AOR05}, the algebra $\mathfrak{L}_5$ is
replaced by a series of algebras. In fact, all the algebras in this series are isomorphic to $\mathfrak{L}_5$ and, by
this reason, they have to be replaced by one rigid algebra.
\item Though there is no algebra that degenerates to $\mathfrak{N}_4$ or $\mathfrak{N}_5$, these algebras belong to 
$\overline{\{O(\mathfrak{N}_3(\alpha))\}_{\alpha\in\mathbb{C}}}$
and hence they are not rigid even in the variety of three-step nilpotent algebras. Thus, they have to
be excluded from the list of rigid algebras obtained in \cite{AOR05}.

\end{enumerate}

For the right classification of nilpotent four-dimensional Leibniz algebras see \cite{alb06} with one excepion of the algebra $\mathfrak{N}_0$ that is missed in this classification and can be found, for example, in \cite{demir}.

\section{Definitions and notation}

All spaces in this paper are considered over $\mathbb{C}$, and we write simply $dim$, $Hom$ and $\otimes$ instead of $dim_{\mathbb{C}}$, $Hom_{\mathbb{C}}$ and $\otimes_{\mathbb{C}}$. An algebra $A$ is a set with a structure of a vector space and a binary operation that induces a bilinear map from $A\times A$ to $A$.

Given an $n$-dimensional vector space $V$, the set $Hom(V \otimes V,V) \cong V^* \otimes V^* \otimes V$ 
is a vector space of dimension $n^3$. This space has a structure of the affine variety $\mathbb{C}^{n^3}.$ Indeed, let us fix a basis $e_1,\dots,e_n$ of $V$. Then any $\mu\in Hom(V \otimes V,V)$ is determined by $n^3$ structure constants $c_{i,j}^k\in\mathbb{C}$ such that
$\mu(e_i\otimes e_j)=\sum\limits_{k=1}^nc_{i,j}^ke_k$. A subset of $Hom(V \otimes V,V)$ is {\it Zariski-closed} if it can be defined by a set of polynomial equations in the variables $c_{i,j}^k$ ($1\le i,j,k\le n$).

Let $T$ be a set of polynomial identities.
All algebra structures on $V$ satisfying polynomial identities from $T$ form a Zariski-closed subset of the variety $Hom(V \otimes V,V)$. We denote this subset by $\mathbb{L}(T)$.
The general linear group $GL(V)$ acts on $\mathbb{L}(T)$ by conjugations:
$$ (g * \mu )(x\otimes y) = g\mu(g^{-1}x\otimes g^{-1}y)$$ 
for $x,y\in V$, $\mu\in \mathbb{L}(T)\subset Hom(V \otimes V,V)$ and $g\in GL(V)$.
Thus, $\mathbb{L}(T)$ is decomposed into $GL(V)$-orbits that correspond to the isomorphism classes of algebras. 
Let $O(\mu)$ denote the orbit of $\mu\in\mathbb{L}(T)$ under the action of $GL(V)$ and $\overline{O(\mu)}$ denote the Zariski closure of $O(\mu)$.

Let $A$ and $B$ be two $n$-dimensional algebras satisfying identities from $T$ and $\mu,\lambda \in \mathbb{L}(T)$ represent $A$ and $B$ respectively.
We say that $A$ degenerates to $B$ and write $A\to B$ if $\lambda\in\overline{O(\mu)}$.
Note that in this case we have $\overline{O(\lambda)}\subset\overline{O(\mu)}$. Hence, the definition of a degeneration does not depend on the choice of $\mu$ and $\lambda$. If $A\not\cong B$, then the assertion $A\to B$ is called a {\it proper degeneration}. We write $A\not\to B$ if $\lambda\not\in\overline{O(\mu)}$.

Let $A$ be represented by $\mu\in\mathbb{L}(T)$. Then  $A$ is  {\it rigid} in $\mathbb{L}(T)$ if $O(\mu)$ is an open subset of $\mathbb{L}(T)$.
 Recall that a subset of a variety is called irreducible if it can't be represented as a union of two non-trivial closed subsets. A maximal irreducible closed subset of a variety is called {\it irreducible component}.
It is well known that any affine variety can be represented as a finite union of its irreducible components in a unique way. We denote by $Rig(\mathbb{L}(T))$ the set of rigid algebras in $\mathbb{L}(T)$.

In particular, $A$ is rigid in $\mathbb{L}(T)$ iff $\overline{O(\mu)}$ is an irreducible component of $\mathbb{L}(T)$.  This is a general fact about algebraic varieties which is well known and was used in many papers (see, for example, \cite{AOR05,B05,kpv}).
For the convenience of the reader we give a proof of this fact here. Let $X$ be an irreducible component of the variety $\mathbb{L}(T)$ such that $O(\mu)\cap X \not= \varnothing$. Then $O(\mu) \subset X$ and hence $\overline{O(\mu)} \subset X$. If $O(\mu)$ is open, then $X = (X \setminus O(\mu)) \cap \overline{O(\mu)}$ and hence $X \subset \overline{O(\mu)}$. Thus, $\overline{O(\mu)} = X$ is an irreducible component. Suppose now that $\overline{O(\mu)}$ is irreducible. Then there exists a closed set $Y \subset \mathbb{L}(T)$ such that $Y \not= \mathbb{L}(T)$ and $\overline{O(\mu)}\cup Y = \mathbb{L}(T)$. In particular, $\overline{O(\mu)} \not\subset Y$ and it follows from the argument above that $Y \cap O(\mu) = \varnothing$. It is well known that $O(\mu)$ is open in $\overline{O(\mu)}$, i.e. there exists an open set
$Z \subset \mathbb{L}(T)$ such that $Z\cap \overline{O(\mu)} = O(\mu)$. Then $O(\mu) = (\mathbb{L}(T)\setminus Y)\cap Z$ is open.
 
An algebra $A$ is called a {\it  Leibniz  algebra}  if it satisfies the identity
$$(xy)z=(xz)y+x(yz).$$ 

An algebra $A$ is called a {\it Zinbiel algebra} if it satisfies the identity 
$$(xy)z=x(yz+zy).$$
It is easy to see that any Lie algebra is a Leibniz algebra.
The classification of nilpotent four-dimensional Leibniz algebras is given in \cite{alb06}.
The classification of four-dimensional Zinbiel algebras is given in \cite{AKO10}.

We denote by $\mathfrak{N}$ the variety of three-step nilpotent associative four-dimensional (Leibniz-Zinbiel) algebras,
by $\mathfrak{Z}$ the variety of four dimensional Zinbiel algebras, and 
by $\mathfrak{L}$ the variety of  nilpotent four dimensional Leibniz algebras. 


We use the following notation: 

\begin{enumerate}
\item $Ann_L(A)=\{ a \in A \mid xa =0 \mbox{ for all } x\in A \}$ is the left  annihilator of $A;$
\item $Ann_R(A)=\{ a \in A \mid ax =0 \mbox{ for all } x\in A \}$ is the right  annihilator of $A;$
\item $Ann(A)=Ann_R(A)\cap Ann_L(A)$ is the  annihilator of $A;$
\item $AZ(A)=\{ a \in A \mid xa +ax=0 \mbox{ for all } x\in A \}$ is the anticommutative center of $A;$%
\item $msub_0(A)$ is a trivial subalgebra of $A$ of the maximal dimension (we fix one for each algebra $A$);
\item $A^{(+2)}$ is the space $\{xy+yx\mid x,y\in A\}$;
\item  $Z(A)$ is the center of $A.$
\end{enumerate}
Given spaces $U$ and $W$, we write simply $U>W$ instead of $dim\,U>dim\,W$. We use also the notation $U\circ W=UW+WU$.


\section{Methods} 

In the present work we use the methods applied to Lie algebras in \cite{BC99,GRH,GRH2,S90}.
First of all, if $A\to B$ and $A\not\cong B$, then $Der(A)<Der(B)$, where $Der(A)$ is the Lie algebra of derivations of $A$. We will compute the dimensions of algebras of derivations and will check the assertion $A\to B$ only for such $A$ and $B$ that $Der(A)<Der(B)$. Secondly, if $A\to C$ and $C\to B$ then $A\to B$. If there is no $C$ such that $A\to C$ and $C\to B$ are proper degenerations, then the assertion $A\to B$ is called a {\it primary degeneration}. If $Der(A)<Der(B)$ and there are no $C$ and $D$ such that $C\to A$, $B\to D$, $C\not\to D$ and one of the assertions $C\to A$ and $B\to D$ is a proper degeneration,  then the assertion $A \not\to B$ is called a {\it primary non-degeneration}. It suffices to prove only primary degenerations and non-degenerations to describe degenerations in the variety under consideration. It is easy to see that any algebra degenerates to the algebra with zero multiplication. From now on we use this fact without mentioning it.


To prove primary degenerations, we will construct families of matrices parametrized by $t$. Namely, let $A$ and $B$ be two algebras represented by the structures $\mu$ and $\lambda$ from $\mathbb{L}(T)$ respectively. Let $e_1,\dots, e_n$ be a basis of $V$ and $c_{i,j}^k$ ($1\le i,j,k\le n$) be the structure constants of $\lambda$ in this basis. If there exist $a_i^j(t)\in\mathbb{C}$ ($1\le i,j\le n$, $t\in\mathbb{C}^*$) such that $E_i^t=\sum\limits_{j=1}^na_i^j(t)e_j$ ($1\le i\le n$) form a basis of $V$ for any $t\in\mathbb{C}^*$, and the structure constants of $\mu$ in the basis $E_1^t,\dots, E_n^t$ are such polynomials $c_{i,j}^k(t)\in\mathbb{C}[t]$ that $c_{i,j}^k(0)=c_{i,j}^k$, then $A\to B$. In this case  $E_1^t,\dots, E_n^t$ is called a {\it parametrized basis} for $A\to B$.

Note also the following fact. Let $B(\alpha)$ be a series of algebras parametrized by $\alpha\in\mathbb{C}$ and $e_1,\dots,e_n$ be a basis of $V$. Suppose also that, for any $\alpha\in\mathbb{C}$, the algebra $B(\alpha)$ can be represented by a structure $\mu(\alpha)\in\mathbb{L}(T)$ having structure constants $c_{i,j}^k(\alpha)\in\mathbb{C}$ in the basis $e_1,\dots,e_n$, where $c_{i,j}^k(t)\in\mathbb{C}[t]$ for all $1\le i,j,k\le n$. Let $A$ be an algebra such that $A\to B(\alpha)$ for $\alpha\in\mathbb{C}\setminus S$, where $S$ is a finite subset of $\mathbb{C}$. Then $A\to B(\alpha)$ for all $\alpha\in\mathbb{C}$. Indeed, if $\lambda\in \mathbb{L}(T)$ represents $A$, then we have $\mu(\alpha)\in\overline{\{\mu(\beta)\}_{\beta\in\mathbb{C}\setminus S}}\subset \overline{O(\lambda)}$ for any $\alpha\in\mathbb{C}$. Thus, to prove that $A\to B(\alpha)$ for all $\alpha\in\mathbb{C}$ we will construct degenerations that are valid for all but finitely many $\alpha$.

Let us describe the methods for proving primary non-degenerations. The main tool for this is the following lemma.

\begin{Lem}[\cite{B99,GRH}]\label{main}
Let $\mathcal{B}$ be a Borel subgroup of $GL(V)$ and $\mathcal{R}\subset \mathbb{L}(T)$ be a $\mathcal{B}$-stable closed subset.
If $A \to B$ and $A$ can be represented by $\mu\in\mathcal{R}$ then there is $\lambda\in \mathcal{R}$ that represents $B$.
\end{Lem}

Since any Borel subgroup of $GL_n(\mathbb{C})$ is conjugated with the subgroup of upper triangular matrices, Lemma \ref{main} can be applied in the following way. Let $A$ and $B$ be two algebras and let $\mu,\lambda$ be the structures in $\mathbb{L}(T)$ representing $A$ and $B$ respectively. Let $Q$ be a set of polynomial equations in variables $x_{i,j}^k$ ($1\le i,j,k\le n$).
Suppose that if $x_{i,j}^k=c_{i,j}^k$ ($1\le i,j,k\le n$) is a solution to all equations in $Q$, then $x_{i,j}^k=\tilde c_{i,j}^k$ ($1\le i,j,k\le n$) is a solution to all equations in $Q$ too in the following cases:
\begin{enumerate}
    \item $\tilde c_{i,j}^k=\frac{\alpha_i\alpha_j}{\alpha_k}c_{i,j}^k$ for some $\alpha_i\in\mathbb{C}^*$ ($1\le i\le n$);
    \item there are some numbers $1\le u<v\le n$ and some $\alpha\in\mathbb{C}$ such that
    $$
    \tilde c_{i,j}^k=\begin{cases}
    c_{i,j}^k,&\mbox{ if $i,j\not=u$ and $k\not=v$},\\
    c_{u,j}^k+\alpha c_{v,j}^k,&\mbox{ if $i=u$, $j\not=u$ and $k\not=v$},\\
    c_{i,u}^k+\alpha c_{i,v}^k,&\mbox{ if $i\not=u$, $j=u$ and $k\not=v$},\\
    c_{i,j}^v-\alpha c_{i,j}^u,&\mbox{ if $i,j\not=u$ and $k=v$},\\
    c_{u,u}^k+\alpha (c_{v,u}^k+c_{u,v}^k)+\alpha^2c_{v,v}^k,&\mbox{ if $i=j=u$ and $k\not=v$},\\
    c_{u,j}^v+\alpha (c_{v,j}^v-c_{u,j}^u)-\alpha^2c_{v,j}^u,&\mbox{ if $i=u$, $j\not=u$ and $k=v$},\\
    c_{i,u}^v+\alpha (c_{i,v}^v-c_{i,u}^u)-\alpha^2c_{i,v}^u,&\mbox{ if $i\not=u$, $j=u$ and $k=v$},\\
    c_{u,u}^v+\alpha (c_{v,u}^v+c_{u,v}^v-c_{u,u}^u)+\alpha^2(c_{v,v}^v-c_{v,u}^u-c_{u,v}^u)-\alpha^3c_{v,v}^u,&\mbox{ if $i=j=u$ and $k=v$}.
    \end{cases}
    $$
\end{enumerate}
Assume that there is a basis $f_1,\dots,f_n$ of $V$ such that the structure constants of $\mu$ in this basis form a solution to all equations in $Q$.
If the structure constants of $\lambda$ don't form a solution to all equations in $Q$ in any basis, then $A\not\to B$.

In particular, it follows from Lemma \ref{main} that $A\not\to B$ in the following cases:
\begin{enumerate}
\item $Ann_L(A)>Ann_L(B)$;
\item $Ann_R(A)>Ann_R(B)$;
\item $Ann(A)>Ann(B)$;
\item $AZ(A)>AZ(B)$;
\item $msub_0(A)>msub_0(B)$;
\item $A^{2}<B^{2}$;
\item $A^{(+2)}<B^{(+2)}$;
\item  $Z(A)>Z(B)$.
\end{enumerate}
In the cases where all of these criteria can't be applied to prove $A\not\to B$, we will define $\mathcal{R}$ by a set of polynomial equations and will give a basis of $V$, in which the structure constants of $\mu$ give a solution to all these equations. We will omit everywhere the verification of the fact that $\mathcal{R}$ is stable under the action of the subgroup of upper triangular matrices and of the fact that $\lambda\not\in\mathcal{R}$ for any choice of a basis of $V$. These verifications can be done by direct calculations.

If the number of orbits under the action of $GL(V)$ on  $\mathbb{L}(T)$ is finite, then the graph of primary degenerations gives the whole picture. In particular, the description of rigid algebras and irreducible components can be easily obtained. Since the variety $\mathfrak{N}$ contains infinitely many non-isomorphic algebras, we have to make some additional work. Let $A(*):=\{A(\alpha)\}_{\alpha\in I}$ be a set of algebras, and let $B$ be another algebra. Suppose that, for $\alpha\in I$, $A(\alpha)$ is represented by the structure $\mu(\alpha)\in\mathbb{L}(T)$ and $B\in\mathbb{L}(T)$ is represented by the structure $\lambda$. Then $A(*)\to B$ means $\lambda\in\overline{\{O(\mu(\alpha))\}_{\alpha\in I}}$, and $A(*)\not\to B$ means $\lambda\not\in\overline{\{O(\mu(\alpha))\}_{\alpha\in I}}$.

Let $A(*)$, $B$, $\mu(\alpha)$ ($\alpha\in I$) and $\lambda$ be as above. To prove $A(*)\to B$ it is enough to construct a family of pairs $(f(t), g(t))$ parametrized by $t\in\mathbb{C}^*$, where $f(t)\in I$ and $g(t)\in GL(V)$. Namely, let $e_1,\dots, e_n$ be a basis of $V$ and $c_{i,j}^k$ ($1\le i,j,k\le n$) be the structure constants of $\lambda$ in this basis. If we construct $a_i^j:\mathbb{C}^*\to \mathbb{C}$ ($1\le i,j\le n$) and $f: \mathbb{C}^* \to I$ such that $E_i^t=\sum\limits_{j=1}^na_i^j(t)e_j$ ($1\le i\le n$) form a basis of $V$ for any  $t\in\mathbb{C}^*$, and the structure constants of $\mu_{f(t)}$ in the basis $E_1^t,\dots, E_n^t$ are such polynomials $c_{i,j}^k(t)\in\mathbb{C}[t]$ that $c_{i,j}^k(0)=c_{i,j}^k$, then $A(*)\to B$. In this case  $E_1^t,\dots, E_n^t$ and $f(t)$ are called a parametrized basis and a {\it parametrized index} for $A(*)\to B$ respectively.

We now explain how to prove $A(*)\not\to B$. Note that if $dim\,Der(A(\alpha))>dim\,Der(B)$ for all $\alpha\in I$ then $A(*)\not\to B$. One can use also the following generalization of Lemma \ref{main}, whose proof is the same as the proof of Lemma \ref{main}.

\begin{Lem}\label{gmain}
Let $\mathcal{B}$ be a Borel subgroup of $GL(V)$ and $\mathcal{R}\subset \mathbb{L}(T)$ be a $\mathcal{B}$-stable closed subset.
If $A(*) \to B$ and for any $\alpha\in I$ the algebra $A(\alpha)$ can be represented by a structure $\mu(\alpha)\in\mathcal{R}$, then there is $\lambda\in \mathcal{R}$ representing $B$.
\end{Lem}

\section{Classification}

We describe all four-dimensional Zinbiel and nilpotent Leibniz algebras in Table 1 below. We use the classification results of \cite{rom,AKO10,alb06,dzhuma5, demir} for this.
Only the algebras $\mathfrak{N}^2_1$ and $\mathfrak{L}_1$ didn't appear in these classifications. The reason for this is that $\mathfrak{N}^2_1$ is a direct sum of two algebras of dimension two and $\mathfrak{L}_1$ is a Lie algebra, while \cite{alb06} contains the classification of indecomposable non-Lie algebras of dimension $4$. The first column of Table 1 contains the names given by us to the algebras under consideration. We use the letters $\mathfrak{N}$, $\mathfrak{Z}$, and $\mathfrak{L}$ in our notation to indicate algebras from $\mathfrak{N}$, $\mathfrak{Z}\setminus\mathfrak{N}$, and $\mathfrak{L}\setminus\mathfrak{N}$ respectively. The second, the third, and the fourth columns of Table 1 contain the names from \cite{rom,AKO10,alb06,dzhuma5}, the multiplication tables, and the dimensions of algebras of derivations respectively for the algebras under consideration. We omit in  multiplication tables zero products of basic elements.

\footnotesize

\begin{center}
Table 1. {\it Zinbiel and nilpotent Leibniz algebras of dimension $4$.}
\begin{equation*}
\begin{array}{|l|l|c|c|} 
\hline
\mbox{notation}  & \mbox{\cite{rom,AKO10,alb06,dzhuma5}}  & \mbox{ multiplication tables } & \mbox{derivations} \\ 

\hline \hline 

\mathfrak{Z}_1^{\mathbb{C}}& W(3)\oplus\mathbb{C} &\begin{array}{c}e_1e_1 = e_2,  e_1e_2 =\frac{1}{2} e_3, e_2e_1=e_3  \end{array} & 6\\ 
\hline
\hline

\mathfrak{Z}_1 & A_1& \begin{array}{c}e_1e_1 = e_2, e_1e_2 = e_3,  e_1e_3 = e_4, e_2e_1 = 2e_3, e_2e_2 = 3e_4, e_3e_1 = 3e_4\end{array}& 4\\ 

\hline
\mathfrak{Z}_2 & A_2& \begin{array}{c}e_1e_1 = e_3, e_1e_2 = e_4,   e_1e_3 = e_4, 
e_3e_1 = 2e_4\end{array}& 5\\ 

\hline
\mathfrak{Z}_3& A_3& \begin{array}{c}e_1e_1 = e_3, e_1e_3 = e_4,  e_2e_2 = e_4, e_3e_1 = 2e_4\end{array} & 4\\ 

\hline
\mathfrak{Z}_4& A_4& \begin{array}{c} e_1e_2 = e_3,  e_1e_3 = e_4, e_2e_1 = -e_3 \end{array}& 5\\ 

\hline
\mathfrak{Z}_5& A_5 & \begin{array}{c}e_1e_2 = e_3, e_1e_3 = e_4, e_2e_1 = -e_3, e_2e_2 = e_4\end{array}& 4 \\ 
\hline
\hline


\mathfrak{N}_1^{\mathbb{C}^2}& L_2\oplus\mathbb{C} & \begin{array}{c}e_1e_1 = e_2 \end{array} 
& 10 \\

\hline
\mathfrak{N}_1^2 & & \begin{array}{c}e_1e_1 = e_2, e_3e_3=e_4  \end{array}& 6 \\ 
\hline

\mathfrak{N}_1^{\mathbb{C}} &L_3\oplus\mathbb{C} & \begin{array}{c}  e_1e_2=  e_3, e_2e_1=-e_3  \end{array}
& 10 \\

\hline
\mathfrak{N}_2^{\mathbb{C}}(\beta)  &L_4(\beta)\oplus\mathbb{C} & \begin{array}{c}   e_1e_1=  e_3, e_1e_2=e_3, e_2e_2=\beta e_3  \end{array}
& 8 \\

\hline
\mathfrak{N}_3^{\mathbb{C}}    &L_5\oplus\mathbb{C}  & \begin{array}{c} \begin{array}{c} e_1e_1=  e_3, e_1e_2=e_3,  e_2e_1=e_3 \end{array}  \end{array}
& 8 \\

\hline
\hline

\mathfrak{N}_0 & & e_1e_2 = e_4, e_3e_1 = e_4& 6   \\ \hline

\mathfrak{N}_1 & R_{12} & \begin{array}{c} e_1e_2 = e_3, e_2e_1 = e_4,  e_2e_2 = -e_3 \end{array} & 5\\

\hline
\mathfrak{N}_2(\gamma )
&R_{13}(\gamma) & \begin{array}{c} e_1e_1 = e_3, e_1e_2 = e_4,  e_2e_1 = -\gamma e_3, e_2e_2 = -e_4 \end{array} & \begin{cases}5&\mbox{ if $\gamma\not=1,$}\\
7&\mbox{ if $\gamma=1;$}
\end{cases}
\\ 
\hline

\mathfrak{N}_3(\alpha    ) & R_{14}(\alpha) & \begin{array}{c}e_1e_1 = e_4, e_1e_2 = \alpha e_4,  e_2e_1 = -\alpha e_4, e_2e_2 = e_4,  e_3e_3 = e_4  \end{array} & \begin{cases}5&\mbox{ if $\alpha\not=0,$}\\
7&\mbox{ if $\alpha=0;$}\\
\end{cases}\\

\hline
\mathfrak{N}_4 &R_{15} & \begin{array}{c} e_1e_2 = e_4, e_1e_3 = e_4, e_2e_1 = -e_4, e_2e_2 = e_4, e_3e_1 = e_4\end{array} & 5 \\ 
\hline
\mathfrak{N}_5 & R_{16} & \begin{array}{c} e_1e_1 = e_4, e_1e_2 = e_4, e_2e_1 = -e_4, e_3e_3 = e_4\end{array}& 5 \\
\hline

\mathfrak{N}_6 & R_{17} & e_1e_2 = e_3, e_2e_1 = e_4 & 6 \\ 
\hline
\mathfrak{N}_7 & R_{11} &  \begin{array}{c} e_1e_1 = e_4, e_1e_2 = e_3, e_2e_1 = -e_3, e_2e_2=2e_3+e_4 \end{array}& 
5 \\
\hline

\mathfrak{N}_9(\alpha) &  R_{20}(\frac{\alpha-1}{\alpha+1}),R_{18}
 & \begin{array}{c} e_1e_2 = e_4, e_2e_1 =\alpha e_4, e_2e_2 = e_3  \end{array} &
7 \\

\hline
\mathfrak{N}_{10} & R_{21} & \begin{array}{c}  e_1e_2 = e_4, e_2e_1 = -e_4, e_3e_3 = e_4 \end{array}& 7 \\

\hline \hline 

\mathfrak{L}_1^{\mathbb{C}} &L_6\oplus\mathbb{C} & \begin{array}{c} e_1e_1 = e_2, e_2e_1=e_3 \end{array} 
& 6 \\ 

\hline
\hline
\mathfrak{L}_1 &  & \begin{array}{c} e_1e_2 = e_3,  e_1e_3 = e_4, e_2e_1 = -e_3, e_3e_1 = -e_4 \end{array}
& 7 \\ 

\hline
\mathfrak{L}_2 & R_1 & \begin{array}{c}e_1e_1 = e_2, e_2e_1 = e_3,  e_3e_1 = e_4 \end{array}
& 4 \\ 

\hline
\mathfrak{L}_3 & R_2& \begin{array}{c}e_1e_1 = e_3, e_1e_2 = e_4, e_2e_1 = e_3,   e_3e_1 = e_4\end{array} 
& 4 \\ 

\hline
\mathfrak{L}_4 & R_3& \begin{array}{c} e_1e_1 = e_3,  e_2e_1 = e_3, e_3e_1 = e_4 \end{array}
& 5 \\ 

\hline
\mathfrak{L}_5  &R_4(0) & \begin{array}{c} e_1e_1 = e_3,  e_2e_1 = e_3, e_2e_2 = e_4, e_3e_1=e_4 \end{array}
& 3 \\ 

\hline
\mathfrak{L}_6 &R_4(1) & \begin{array}{c} e_1e_1 = e_3, e_1e_2 =  e_4, e_2e_1 = e_3, e_2e_2 = e_4, e_3e_1=e_4 \end{array}
& 4 \\

\hline
\mathfrak{L}_7 &R_5 & \begin{array}{c} e_1e_1 = e_3, e_1e_2 = e_4,  e_3e_1 = e_4  \end{array}
& 5 \\ 

\hline
\mathfrak{L}_8  &R_6 & \begin{array}{c} e_1e_1 = e_3,   e_2e_2 = e_4, e_3e_1=e_4  \end{array}
& 4 \\

\hline
\mathfrak{L}_9  &R_8  & \begin{array}{c} e_1e_2=-e_3+e_4, e_1e_3=-e_4,  e_2e_1 = e_3,  e_3e_1=e_4   \end{array}
& 5 \\

\hline
\mathfrak{L}_{10} &R_9 & \begin{array}{c} e_1e_2=-e_3, e_1e_3=-e_4, e_2e_1 = e_3,  e_2e_2=e_4,  e_3e_1=e_4   \end{array}
& 5  \\ 

\hline
\mathfrak{L}_{11}  & R_{10} & \begin{array}{c} e_1e_1 = e_4,  e_1e_2=-e_3, e_1e_3=-e_4,  e_2e_1 = e_3, e_2e_2=e_4,  e_3e_1=e_4   \end{array}
& 4 \\ 

\hline
\mathfrak{L}_{12} &R_7 & \begin{array}{c} e_1e_1 = e_4,  e_1e_2=-e_3, e_1e_3=-e_4, e_2e_1 = e_3,   e_3e_1=e_4   \end{array}
& 6  \\

\hline
\end{array}
\end{equation*}
\end{center}

\normalsize 

Since $\mathfrak{N}_3(\alpha)\cong\mathfrak{N}_3(-\alpha)$ for any $\alpha\in\mathbb{C}$, we assume everywhere that $\alpha$ in the notation $\mathfrak{N}_3(\alpha)$ satisfies either the condition $Re(\alpha)>0$ or the conditions $Re(\alpha)=0$ and $Im(\alpha)\ge 0$. With this assumption, different structures from Table 1 define non-isomorphic algebras.

\section{Degenerations}

\begin{Th}\label{third}\label{theorem}
The graph of primary degenerations for Zinbiel and nilpotent Leibniz algebras of dimension 4 has the following form: 
\end{Th}
 
 \scriptsize
 
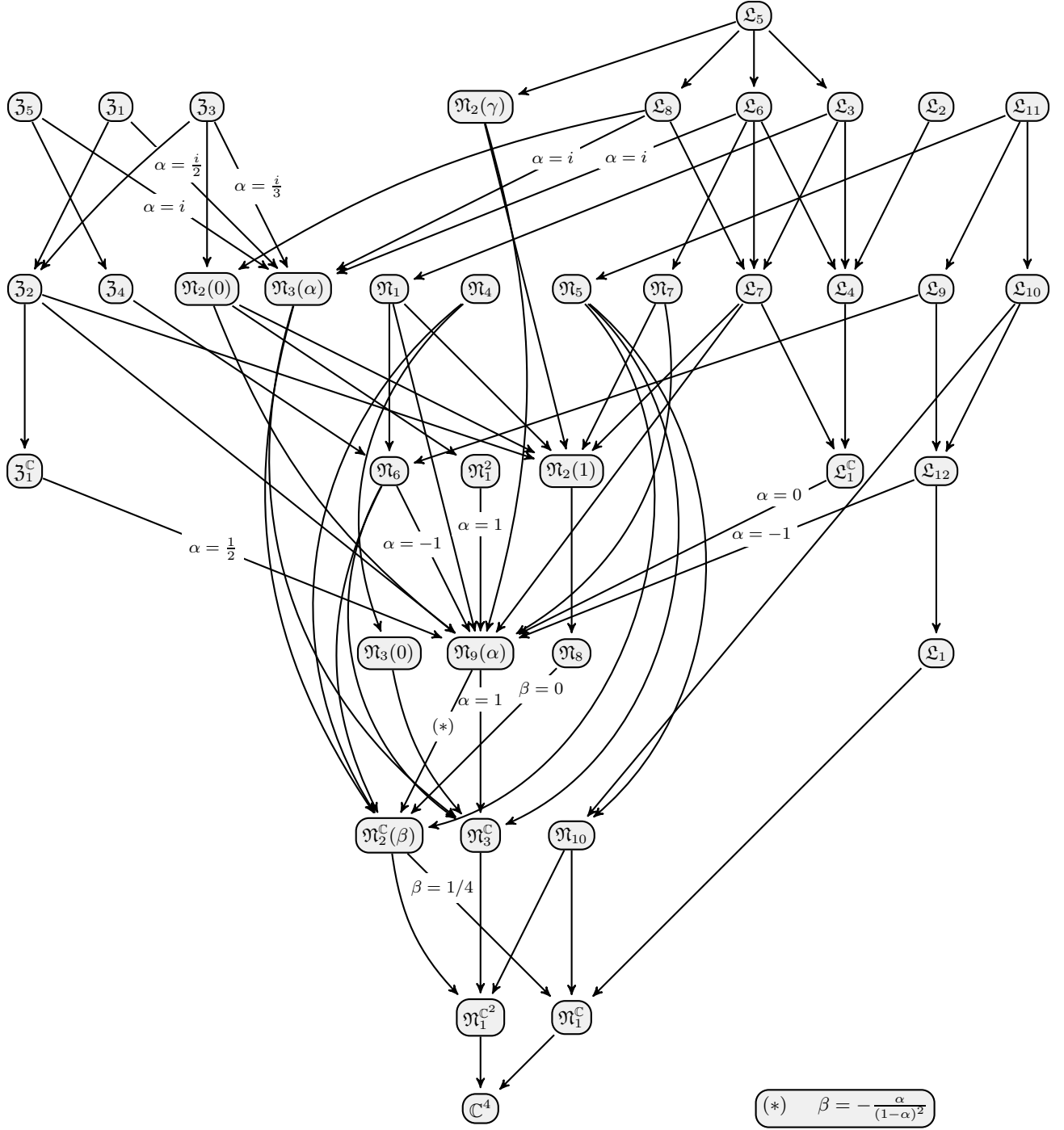
\begin{center}
\begin{tikzpicture}[->,>=stealth',shorten >=0.08cm,auto,node distance=1.5cm,
                    thick,main node/.style={rectangle,draw,fill=gray!12,rounded corners=1.5ex,font=\sffamily \bf 
                    \bfseries },rigid node/.style={rectangle,draw,fill=black!20,rounded corners=1.5ex,font=\sffamily \tiny \bfseries },style={draw,font=\sffamily \scriptsize \bfseries }]
                    
\node (103)   {};

\node (1031)[right of=103]{};    
\node (10311)[right of=103]{};    
\node (103111)[right of=10311]{};    
\node (1031111)[right of=103111]{};    
\node (10311111)[right of=1031111]{};    
\node (103111111)[right of=10311111]{};    
\node (1031111111)[right of=103111111]{};    
\node (10311111111)[right of=1031111111]{};    
\node (103111111111)[right of=10311111111]{};    

\node  [main node] (l5) [right of=10311111] {$\mathfrak{L}_5$};

\node (104) [below          of=103]       {};
\node (1005) [below          of=104]      {};
\node (105) [below          of=1005]      {};
\node (1006) [below          of=105]      {};
\node (106) [ below         of=1006]      {};
\node (1006) [below          of=106]      {};
\node (107) [ below         of=1006]      {};
\node (1007) [below          of=107]      {};
\node (108) [ below         of=1007]      {};
\node (1008) [below          of=108]      {};
\node (110) [ below         of=1008]      {};
\node (116) [ below         of=110]       {};

\node [main node] (z5)   [below of =103]                        {$\mathfrak{Z}_5$ };
\node [main node] (z1)   [right of =z5]                        {$\mathfrak{Z}_1$ };

\node [main node] (z3)   [right of =z1]                        {$\mathfrak{Z}_3$ };

\node (1041) [right of=z3]{};    
\node (1042) [right of=1041]{};

\node [main node] (l8)  [right of =1042  ]                      {$\mathfrak{L}_8$ };
\node (1043) [right of=l8]{};    
\node [main node] (l6)  [right of =1043  ]                       {$\mathfrak{L}_6$ };
\node [main node] (l3)   [right of =l6  ]                       {$\mathfrak{L}_3$ };
\node [main node] (l2)   [right of =l3]                          {$\mathfrak{L}_2$ };

\node [main node] (l11)  [right of =l2  ]                      {$\mathfrak{L}_{11}$ };

\node [main node] (z2)   [below of =1005]                        {$\mathfrak{Z}_2$ };
\node [main node] (z4)   [right of =z2]                        {$\mathfrak{Z}_4$ };

\node [main node] (n3a)  [right of =z4]                      {$\mathfrak{N}_3(\alpha)$ };

\node [main node] (n20)  [right of =n3a]                       {$\mathfrak{N}_2(\gamma)$ };
\node [main node] (n1)   [right of =n20]                        {$\mathfrak{N}_1$ };
\node [main node] (n4)   [right of =n1]                      {$\mathfrak{N}_4$ };
\node [main node] (n5)   [right of =n4]                        {$\mathfrak{N}_5$ };
\node [main node] (n7)   [right of =n5]                       {$\mathfrak{N}_7$ };

\node (1052)[right of=n5]{};

\node [main node] (l7)   [right of =n7]                          {$\mathfrak{L}_7$ };
\node [main node] (l4)   [right of =l7]                        {$\mathfrak{L}_{4}$ };
\node [main node] (l9)   [right of =l4]                          {$\mathfrak{L}_9$ };
\node [main node] (l10)  [right of =l9  ]                      {$\mathfrak{L}_{10}$ };


\node (1061)[right of=106]{};    
\node (10611)[right of=106]{};    
\node (106111)[right of=10611]{};

\node  [main node](z1c) [left of=1061] {$\mathfrak{Z}_1^{\mathbb{C}}$};    

\node [main node] (n6)   [right of =106111]                        {$\mathfrak{N}_6$ };

\node (1061111)[right of=n6]{};    

\node [main node] (n0)   [right of =1061111]                        {$\mathfrak{N}_0$  };

\node (1061111x) [right of=n0]{};

\node [main node] (n12)   [right of =1061111x]                       {$\mathfrak{N}_1^2$ };

\node (1062) [right of=n12]{};

\node [main node] (l1c)   [right of =1062]                        {$\mathfrak{L}_1^{\mathbb{C}}$ };
\node (1063)[right of=l1c]{};    
 
\node [main node] (l12)   [right of =1063]                        {$\mathfrak{L}_{12}$ };

\node (1071)[right of=107]{};    
\node (10711)[right of=107]{};    
\node (107111)[right of=10711]{};

\node [main node] (n30)   [right of =107111]                        {$\mathfrak{N}_3(0)$ };

\node (1071111)[right of=n30]{};    
\node [main node] (n9a)   [right of =1071111]                        {$\mathfrak{N}_9(\alpha)$ };
\node [main node] (n21)   [right of =n9a]                         {$\mathfrak{N}_2(1)$ };
\node (12xx)[right of=n21]{};    

\node [main node] (n10)   [right of =12xx]                        {$\mathfrak{N}_{10}$ };

\node (1071)[right of=n21]{};    
\node (1072)[right of=1071]{};    
\node (1073)[right of=1072]{};

\node [main node] (l1)   [right of =1073]                          {$\mathfrak{L}_{1}$ };

\node (1081)[right of=108]{};    
\node (10811)[right of=108]{};    
\node (108111)[right of=10811]{};    
\node (1081111)[right of=108111]{};

\node [main node] (n3c)  [right of =1081111]                            {$\mathfrak{N}_{3}^{\mathbb{C}}$ };
\node (10xxx)[right of=n3c]{};    
\node (10xxxx)[right of=10xxx]{};    

\node [main node] (n2cb)  [right of =10xxxx]                        {$\mathfrak{N}_{2}^{\mathbb{C}}(\beta)$ };

\node (1011)[right of=110]{};    

\node (1021)[right of=1011]{};    
\node (10211)[right of=1021]{};    
\node (102111)[right of=1021]{};    
\node (1021111)[right of=102111]{};

\node [main node] (n1c2)   [right of =1021111]                        {$\mathfrak{N}_{1}^{\mathbb{C}^2}$ };
\node [main node] (n1c)   [right of =n1c2]                           {$\mathfrak{N}_{1}^{\mathbb{C}}$ };

\node (101)[right of=116]{};    
\node (1011)[right of=101]{};    
\node (10111)[right of=1011]{};    
\node (101111)[right of=1011]{};    
\node (1011111)[right of=101111]{};    

    \node[main node] (c4) [ right  of=1011111]       {$\mathbb{C}^4$};
 
\node (191)[right of=c4]{};    
\node (1911)[right of=191]{};    
\node (19111)[right of=1911]{};    
\node (191111)[right of=19111]{};    
\node (1911111)[right of=191111]{};    
 
    \node[main node] (ll) [ right  of=19111]       {$(*)\hspace{0.5cm}\beta=-\frac{\alpha}{(1-\alpha)^2}$};

\path[every node/.style={font=\sffamily\small}]

(z1)  edge   node[above=-7, left=-20,  fill=white]{\tiny $\alpha=\frac{i}{2}$}   node{}  (n3a) 

(z1)  edge   node[above=9, left=-8,  fill=white]{\tiny $\gamma=\frac{1}{9}$}   node{}  (n20)

(z1)  edge  node{}   (z2)
 
(z2)  edge  node{}   (z1c)

(z2)  edge  [bend right=-13] node{}   (n21)

(z2)  edge  [bend left=-12] node{}   (n9a)
(z2)  edge node{}  (n0)

(z3)  edge  [bend left=30] node[above=9, left=-20, fill=white]{\tiny  $\alpha= \frac{i}{3}$}   node{}  (n3a) 
(z3)  edge  [bend right=0] node{}   (z2)
(z3)  edge  [bend left=30] node[above=20, left=-5, fill=white]{\tiny  $\gamma= 0$}   node{}    (n20)

(z4)  edge node{}  (n6) 
(z4)  edge node{}  (n0)

(z5)  edge [bend left=-0]   node{}  (z4) 

(z5)  edge  [bend left=-15] node[below=16, right=3, fill=white]{\tiny $\alpha=  i$}   node{}  (n3a)

(n1)  edge node{}  (n6) 
(n1)  edge [bend right= 11] node{}  (n21)
(n1)  edge node{}  (n9a)

(z1c)  edge  node[above=5, left, fill=white]{\tiny  $\alpha=\frac{1}{2}$}   node{}  (n9a)

(n20)  edge node{}  (n12) 

(n20)  edge node{}  (n21)

(n20)  edge [bend right=0] node{}  (n9a) 

(n30)  edge [bend right=0] node{}  (n3c)

(n3a)  edge node{}  (n0)

(n4)  edge [bend right=-25]node{}  (n30) 
(n4)  edge node{}  (n0)

(n5)  edge [bend left=15] node{}  (n10) 
(n5)  edge node{}  (n0)

(n6)  edge node[above=10, left=-6, fill=white]{\tiny $\alpha=-1$} node{}  (n9a) 

(n6)  edge [bend right=27] node{}  (n2cb)

(n6)  edge   node{}  (n3c) 

(n7)  edge  node{}  (n21) 

(n7)  edge  node{}  (n9a)

(n21)  edge [bend left=0] node[above=0, right=-25, fill=white]{\tiny \ \ \ \  $\beta=0$} node{} (n2cb)

(n9a)  edge  [bend left=0]  node[above=22, right=-2, fill=white]{\tiny $\alpha=1$} node{}   (n3c) 

(n9a)  edge   node[above=8, right=-17, fill=white]{\tiny $(*)$} node{}   (n2cb) 

(n10)  edge [bend left=00] node[above=0, left=-25, fill=white]{\tiny $\beta=1/4$}  (n2cb)

(n1c2)  edge node{}  (c4) 
    
(n12)  edge  [bend left= 12] node[above=22,left=-38, fill=white]{\tiny $\alpha=1$} node{}  (n9a) 

(n1c)  edge node{}  (c4)

(n2cb)  edge  node[above=0, left=-17, fill=white]{\tiny $\beta=1/4$}  (n1c) 
(n2cb)  edge [bend right=0] node{}  (n1c2) 
(n3c)   edge node{}  (n1c2)

(l1)  edge node{}  (n1c)

(l2)  edge node{}  (l4) 

(l3)  edge node{}  (l4)

(l3)  edge node{}  (l7)

(l3)  edge node{}  (n1)

(l4)  edge node{}  (l1c)

(l5)  edge node{}  (l3)

(l5)  edge node{}  (l6)

(l5)  edge node{}  (l8) 

(l5)  edge [bend right=30] node{}  (n20)

(l6)  edge  [bend left=0]  node[above=3, right=0, fill=white]{\tiny $\alpha=  i$}   node  {} (n3a)

(l6)  edge node{}  (n7) 

(l6)  edge node{}  (l4) 

(l6)  edge node{}  (l7) 

(l7)  edge node{}  (l1c) 
(l7)  edge [bend right=13] node{}  (n21) 
(l7)  edge node{}  (n9a) 
(l7)  edge node{}  (n0)

(l8)  edge node{}  (l7) 
(l8)  edge [bend left=0] node[above=3, right=-8, fill=white]{\tiny $\gamma= 0$}   node  {}  (n20) 
(l8)  edge [bend left=0] node[above=3, right=-15, fill=white]{\tiny $\alpha= i$}   node  {} (n3a)

(l9)  edge node{}  (n6)
(l9)  edge node{}  (n0)
(l9)  edge node{}  (l12)

(l6)  edge  [bend left=0]  node[above=3, right=0, fill=white]{\tiny $\alpha=  i$}   node  {} (n3a)

(l10)  edge [bend left=00] node{}  (n10)

(l10)  edge node{}  (l12)

(l11)  edge node{}  (n5)

(l11)  edge node{}  (l9)

(l11)  edge node{}  (l10)
 
(l12)  edge node[above=5, right= -10, fill=white]{\tiny $\alpha= -1$} node{}  (n9a)

(l12)  edge node{}  (l1)
 
(l1c)  edge node[above=25, left=-65, fill=white]{\tiny \ \ \ \  $\alpha=0$}   node  {} (n9a)


(n0)  edge node{}  (n3c)
(n0)  edge [bend right=-20] node{}  (n2cb)

        ;
\end{tikzpicture}

\end{center}

\normalsize

 \begin{Proof} We prove all required primary degenerations in Table 2 below. Let us consider the first degeneration $\mathfrak{Z}_1^{\mathbb{C}} \to  \mathfrak{N}_9 \left(\frac{1}{2}\right)$ to clarify this table. Write nonzero products in $\mathfrak{Z}_1^{\mathbb{C}}$ in the basis  $E_i^t$: 
$$E_1^tE_2^t=te_3=E_4^t;\ E_2^tE_2^t=t^2e_2=t^2E_1^t+E_3^t;\  E_2^tE_1^t=\frac{1}{2}te_3=\frac{1}{2}E_4^t.$$
 It is easy to see that for $t=0$ we obtain the multiplication table of $\mathfrak{N}_9 \left(\frac{1}{2}\right)$. The remaining degenerations can be considered in the same way.

\begin{center} \footnotesize Table 2. {\it Primary degenerations of Zinbiel and nilpotent Leibniz algebras of dimension $4$.}
$$\begin{array}{|l|c|}
\hline
\mbox{degenerations}  &  \mbox{parametrized bases}\\
\hline
\hline

\mathfrak{Z}_1^{\mathbb{C}} \to  \mathfrak{N}_9 \left(\frac{1}{2}\right)   &  E_1^t=e_2-e_4, E_2^t=te_1, E_3^t=t^2e_4, E_4^t=te_3  \\
\hline
\hline

\mathfrak{Z}_1 \to  \mathfrak{Z}_2   &  E_1^t=te_1, E_2^t=3t^2e_2-2te_3+2e_4, E_3^t=t^2e_2, E_4^t=t^2e_4  \\
\hline

\mathfrak{Z}_1 \to  \mathfrak{N}_2\left(\frac{1}{9}\right)   &  E_1^t=3te_1+2te_2+2te_3+3te_4, E_2^t=te_1, E_3^t=9t^2e_3+\frac{45t^2}{2}e_4, E_4^t=t^2e_3+\frac{3t^2}{2}e_4  \\
\hline

\mathfrak{Z}_1 \to  \mathfrak{N}_3\left(\frac{i}{2}\right)   &  E_1^t=\frac{t}{2}(e_1+e_3), E_2^t=\frac{it}{2}(e_1-e_3), E_3^t=\frac{t}{\sqrt{3}}e_2, E_4^t=t^2e_4  \\
\hline
\hline

\mathfrak{Z}_2 \to  \mathfrak{Z}_1^{\mathbb{C}}   &  E_1^t=e_1, E_2^t=e_3, E_3^t=2e_4, E_4^t=te_2  \\
\hline

\mathfrak{Z}_2\to \mathfrak{N}_{0}  &  E_1^t=te_1,  E_2^t=e_2,  E_3^t= - \frac{1}{2}e_2+\frac{1}{2} e_3, E_4^t = t e_4\\

\hline

\mathfrak{Z}_2 \to  \mathfrak{N}_2(1)   &  E_1^t=te_1, E_2^t=e_2-te_1, E_3^t=t^2e_3, E_4^t=te_4-t^2e_3  \\
\hline

\mathfrak{Z}_2 \to  \mathfrak{N}_9(\alpha)   &  E_1^t=(2\alpha-1)e_2+e_3, E_2^t=te_1, E_3^t=t^2e_3, E_4^t=2te_4  \\
\hline
\hline

\mathfrak{Z}_3 \to  \mathfrak{Z}_2   &  E_1^t=te_1+2e_2, E_2^t=t^3e_2-t^2e_3, E_3^t=t^2e_3+4e_4, E_4^t=t^3e_4  \\
\hline

\mathfrak{Z}_3 \to  \mathfrak{N}_2(0)   &  E_1^t=t^2e_1-2te_2, E_2^t=te_2+e_3, E_3^t=t^4e_3+4t^2e_4, E_4^t=-t^2e_4  \\
\hline

\mathfrak{Z}_3 \to  \mathfrak{N}_3\left(\frac{i}{3}\right)   & E_1^t=\frac{2t^4}{3}e_1, E_2^t=\frac{2it^4}{3}e_1-it^2e_3+\frac{9i}{4}e_4, E_3^t=t^3e_2, E_4^t=t^6e_4  \\
\hline
\hline

\mathfrak{Z}_4\to \mathfrak{N}_{0}  &  E_1^t=te_1,  E_2^t=e_3-e_4,  E_3^t= - e_2+ e_3, E_4^t = t e_4\\

\hline

\mathfrak{Z}_4 \to  \mathfrak{N}_6   & E_1^t=te_1, E_2^t=e_2+e_3, E_3^t=t(e_3+e_4), E_4^t=-te_3  \\
\hline
\hline

\mathfrak{Z}_5 \to  \mathfrak{Z}_{4}   & E_1^t=e_1, E_2^t=te_2, E_3^t=te_3, E_4^t=te_4  \\
\hline

\mathfrak{Z}_5 \to  \mathfrak{N}_3(i)   & E_1^t=t(e_1+e_3), E_2^t=-it(e_1-e_3), E_3^t=te_2, E_4^t=t^2e_4  \\
\hline
\hline

\mathfrak{N}_1^2  \to  \mathfrak{N}_{9}(1)    & 
E_1^t=te_1, E_2^t=e_1+e_3, E_3^t=e_2+e_4, E_4^t=te_2 \\
\hline
\hline

\mathfrak{N}_{2}^{\mathbb{C}}(\beta) \to  \mathfrak{N}_1^{\mathbb{C}^2}    & E_1^t=e_1, E_2^t=e_3, E_3^t=te_2, E_4^t=e_4 \\
\hline

\mathfrak{N}_{2}^{\mathbb{C}}(\frac{1}{4}) \to  \mathfrak{N}_1^{\mathbb{C}}    &
E_1^t=te_1-2te_2, E_2^t=2t^2e_2, E_3^t=t^3e_3, E_4^t=e_4 \\
\hline
\hline

\mathfrak{N}_{3}^{\mathbb{C}} \to  \mathfrak{N}_1^{\mathbb{C}^2}    & E_1^t=e_1, E_2^t=e_3, E_3^t=te_2, E_4^t=e_4 \\
\hline
\hline

\mathfrak{N}_0 \to  \mathfrak{N}_2^{\mathbb{C}}(\beta)   & 
E_1^t=e_1 +  e_2, E_2^t=ue_1+e_2  - ue_3, E_3^t=e_4, E_4^t=t e_3, \mbox{ where $u-u^2=\beta$}    \\
\hline

\mathfrak{N}_0 \to  \mathfrak{N}_3^{\mathbb{C}}   & E_1^t=e_1 +  e_2, E_2^t=e_2 + e_3, E_3^3=e_4, E_4^t=te_3  \\
\hline
\hline

\mathfrak{N}_1 \to  \mathfrak{N}_2(1)   & 
E_1^t=-ite_2, E_2^t=t^2e_1+ite_2, E_3^t=t^{2}e_3, E_4^t=-it^3e_4-t^2e_3  \\
\hline

\mathfrak{N}_1 \to  \mathfrak{N}_6   & E_1^t=e_1, E_2^t=te_2, E_3^t=te_3, E_4^t=te_4  \\
\hline

\mathfrak{N}_1 \to  \mathfrak{N}_9(\alpha)   & 
E_1^t=\frac{t}{u}e_1, E_2^t=\frac{1}{u}e_1+ue_2, E_3^t=e_4-\alpha e_3, E_4^t=te_3, \mbox{ where $u^2=\alpha+1$}  \\
\hline
\hline

\mathfrak{N}_2(0) \to  \mathfrak{N}_1^2   & 
E_1^t=te_1, E_2^t=t^2e_3, E_3^t=ie_2, E_4^t=e_4 \\
\hline

\mathfrak{N}_2(1) \to  \mathfrak{N}_2^{\mathbb{C}}(0)    & E_1^t=te_2, E_2^t=te_1, E_3^t=t^2e_3, E_4^t=e_4-e_3 \\
\hline

\mathfrak{N}_2(\gamma \neq 1) \to  \mathfrak{N}_2(1)    & 
E_1^t=e_1, E_2^t=-e_1+te_2, E_3^t=e_3, E_4^t=-e_3+te_4 \\
\hline


\mathfrak{N}_2(\gamma \neq 1) \to  \mathfrak{N}_9(\alpha)& 
\begin{array}{c}
E_1^t=(1+\gamma u)te_2, E_2^t=e_1-ue_2, E_3^t=(1+\gamma u)e_3-u(u+1)e_4, E_4^t=(1-\gamma)ute_4,\\
\mbox{ where $(1+u)(1+\gamma u)=(1-\gamma)\alpha u$}
\end{array}\\
\hline
\hline

\mathfrak{N}_3(\alpha\neq 0) \to  \mathfrak{N}_2^{\mathbb{C}}(\beta)    & 
\begin{array}{c}E_1^t=e_1,
E_2^t=\frac{1}{2}e_1+\frac{1}{2\alpha}e_2+\frac{u}{2\alpha} e_3,
E_3^t=e_4, E_4^t=te_3,\mbox{ where $u^2=4\alpha^2\beta-\alpha^2-1$}\\
\end{array} \\
\hline

\mathfrak{N}_3(\alpha ) \to  \mathfrak{N}_3^{\mathbb{C}}    & 
E_1^t=e_1, E_2^t=e_1+ie_3, E_3^t=e_4, E_4^t=te_2 \\
\hline

\mathfrak{N}_3(\alpha\neq 0)\to \mathfrak{N}_{0}  &  E_1^t=e_1 + i e_3,  E_2^t=- \frac{t}{2\alpha}  e_2+\frac{i t}{2} e_3,  E_3^t= \frac{t}{2\alpha} e_2+ \frac{i t}{2} e_3, E_4^t = -t e_4\\
\hline
\hline

\mathfrak{N}_4 \to  \mathfrak{N}_{2}^{\mathbb{C}}(\beta)   & 
\begin{array}{c}E_1^t=\frac{t}{2u}e_1+\frac{t}{4u}e_2+ ut(1-\frac{1}{16\beta})e_3, E_2^t=ut\left(e_2 + \frac{e_3}{2}\right),\\ E_3^t=t^2e_4, E_4^t=t^2e_3,\mbox{ where $u^2=\beta$}\end{array} \\
\hline

\mathfrak{N}_4\to \mathfrak{N}_{0}  &  E_1^t=e_1 + e_2 - \frac{1}{2} e_3,  E_2^t= \frac{t}{2} e_2,  E_3^t= - \frac{t}{2} e_2+t e_3, E_4^t = t e_4\\
\hline

\mathfrak{N}_4 \to  \mathfrak{N}_3(0)   & 
E_1^t=t^2e_1+it e_2+e_3, E_2^t=te_2-ie_3, E_3^t=it^2e_1-te_2, E_4^t=t^{2}e_4 \\
\hline
\hline

\mathfrak{N}_5 \to  \mathfrak{N}_2^{\mathbb{C}}(\beta)   
& E_1^t=t e_1,  
E_2^t=\frac{1}{2} t(e_1+e_2)+\frac{u}{2}t e_3,  E_3^t=t^{2}e_4, E_4^t=t^{2}e_2,\mbox{ where $u^2=4\beta-1$} \\
\hline

\mathfrak{N}_5 \to  \mathfrak{N}_3^{\mathbb{C}}   & 
E_1^t=e_1, E_2^t=e_1+ie_3, E_3^t=e_4, E_4^t=te_2 \\
\hline

\mathfrak{N}_5\to \mathfrak{N}_{0}  &  E_1^t=i e_1  +  e_3,  E_2^t= - \frac{i t}{2} e_2+\frac{t}{2}  e_3,  E_3^t=\frac{i t}{2} e_2+ \frac{t}{2}  e_3, E_4^t = t e_4\\
\hline

\mathfrak{N}_5 \to  \mathfrak{N}_{10}   & E_1^t=t^2e_1, E_2^t=e_2, E_3^t=te_3, E_4^t=t^2e_4 \\
\hline
\hline

\mathfrak{N}_6 \to  \mathfrak{N}_{2}^{\mathbb{C}}\left(\beta\right)   & 
\begin{array}{c}
E_1^t=t(e_1+e_2), E_2^t=t(u e_1+(1-u)e_2), E_3^t=t^{2}(e_3+e_4), E_4^t=u e_3+(1-u)e_4,\\
\mbox{where }u^2-u+\beta=0
\end{array}\\
\hline

\mathfrak{N}_6 \to  \mathfrak{N}_3^{\mathbb{C}}   & E_1^t=t(e_1+e_2), E_2^t=2te_2, E_3^t=t^2(e_3+e_4), E_4^t=e_3-e_4 \\
\hline

\mathfrak{N}_6 \to  \mathfrak{N}_9(-1)   
& E_1^t=t(e_1+e_2), E_2^t=e_1-e_2, E_3^t=-(e_4+e_3), E_4^t=t(e_4-e_3) \\
\hline
\hline


\mathfrak{N}_7 \to  \mathfrak{N}_2(1)   & 
E_1^t=e_1+e_2, E_2^t=-e_1+(t-1)e_2, E_3^t=2e_3+2e_4, E_4^t=2(t-1)e_3-2e_4 \\
\hline

\mathfrak{N}_7 \to  \mathfrak{N}_9(\alpha)   & 
\begin{array}{c}
E_1^t=te_1, E_2^t=\frac{u(\alpha+1)}{2}e_1+e_2, E_3^t=2e_3+u(\alpha-1)e_4, E_4^t=ute_4,\\
\mbox{where }\left(\frac{\alpha+1}{2}\right)^2u^2-(\alpha-1)u+1=0
\end{array}\\
\hline
\hline

\mathfrak{N}_9(\alpha\neq 1) \to  \mathfrak{N}_{2}^{\mathbb{C}}\left(-\frac{\alpha}{(1-\alpha)^2}\right)
& \begin{array}{c} E_1^t=te_2, 
E_2^t=\frac{1}{\alpha-1}e_1+ \frac{t}{1-\alpha} e_2,
E_3^t=te_4, E_4^t=te_3-e_4
\end{array} \\
\hline


\mathfrak{N}_9(1) \to  \mathfrak{N}_{3}^{\mathbb{C}}   & 
E_1^t=t\left(e_1+\frac{e_2}{2}\right), E_2^t=te_2, E_3^t=t^2e_4, E_4^t=e_3 \\
\hline
\hline



\mathfrak{N}_{10} \to  \mathfrak{N}_{2}^{\mathbb{C}}(\frac{1}{4})    & 
E_1^t=2(e_1+e_3), E_2^t=e_2+e_3, E_3^t=4e_4, E_4^t=te_3 \\
\hline
\hline

\mathfrak{L}_{1}^{\mathbb{C}} \to  \mathfrak{N}_9(0)    & 
E_1^t=t^3e_2-te_3+e_4, E_2^t=t^2e_1, E_3^t=t^2e_3-te_4, E_4^t=t^4e_4 \\
\hline
\hline

\mathfrak{L}_{1} \to  \mathfrak{N}_1^{\mathbb{C}}    & E_1^t=e_1, E_2^t=te_2, E_3^t=te_3, E_4^t=e_4 \\
\hline
\hline

\mathfrak{L}_{2} \to  \mathfrak{L}_4     & E_1^t=t(e_1+e_2), E_2^t=te_2, E_3^t=t^2e_3, E_4^t=t^3e_4 \\
\hline
\hline

\mathfrak{L}_{3} \to  \mathfrak{N}_1     & 
E_1^t=-2te_2, E_2^t=t(e_1+e_2-e_3), E_3^t=-2t^2e_3, E_4^t=-2t^2e_4 \\
\hline

\mathfrak{L}_{3} \to  \mathfrak{L}_4     & 
E_1^t=t^{-1}e_1, E_2^t=t^{-1}e_2, E_3^t=t^{-2}e_3, E_4^t=t^{-3}e_4 \\
\hline

\mathfrak{L}_{3} \to  \mathfrak{L}_7     & 
E_1^t=te_1, E_2^t=t^2e_2, E_3^t=t^2e_3, E_4^t=t^{3}e_4 \\
\hline
\hline

\mathfrak{L}_{4} \to  \mathfrak{L}_1^{\mathbb{C}}    & E_1^t=e_1, E_2^t=e_3, E_3^t=e_4, E_4^t=te_2 \\
\hline

\end{array}$$
\end{center}

\begin{center}\footnotesize
$$\begin{array}{|l|c|}
\hline
\mbox{degenerations}  &  \mbox{parametrized bases}\\
\hline
\hline

\mathfrak{L}_5 \to  \mathfrak{N}_2(\gamma)     & 
\begin{array}{c} E_1^t=(\gamma-1)te_2+\gamma te_3,
E_2^t=(\gamma-1)^2te_1-(\gamma-1)\gamma t e_2-\gamma^2 te_3,\\
E_3^t=(\gamma-1)^2t^2e_4,
E_4^t=(\gamma-1)^3t^2e_3\end{array}\\
\hline

\mathfrak{L}_5 \to  \mathfrak{L}_3     & 
E_1^t=-e_1+(t+1)e_2+(2t+1)e_3, E_2^t=t(e_2+e_3), E_3^t=-te_3, E_4^t=te_4 \\
\hline

\mathfrak{L}_5 \to  \mathfrak{L}_6     & 
\begin{array}{c} E_1^t=t^{-1}e_1+t^{-2}e_2,
E_2^t=t^{-2}e_2, E_3^t=t^{-3}e_3+t^{-4}e_4, E_4^t=t^{-4}e_4 
\end{array} \\
\hline

\mathfrak{L}_5 \to  \mathfrak{L}_8     & 
E_1^t=t^2e_1, E_2^t=t^3e_2, E_3^t=t^4e_3, E_4^t=t^6e_4 \\
\hline
\hline

\mathfrak{L}_6 \to  \mathfrak{N}_3(i)     & 
E_1^t=t^4e_1, E_2^t=it^4e_1-2it^2e_3+2ie_4, E_3^t=t^3e_2-te_3, E_4^t=t^8e_3 \\
\hline

\mathfrak{L}_6 \to  \mathfrak{N}_7     & 
E_1^t=te_1, E_2^t=t(-e_1+2e_2-2e_3), E_3^t=t^2(2e_4-e_3), E_4^t=t^{2}e_3 \\
\hline

\mathfrak{L}_6 \to  \mathfrak{L}_4     & 
E_1^t=t^{-1}e_1, E_2^t=t^{-1}e_2, E_3^t=t^{-2}e_3, E_4^t=t^{-3}e_4 \\
\hline

\mathfrak{L}_6 \to  \mathfrak{L}_7     & 
E_1^t=te_1, E_2^t=t^2e_2, E_3^t=t^2e_3, E_4^t=t^{3}e_4 \\
\hline
\hline

\mathfrak{L}_7 \to  \mathfrak{N}_0     & 
E_1^t=te_1, E_2^t=e_2, E_3^t=e_3, E_4^t=te_4 \\
\hline

\mathfrak{L}_7 \to  \mathfrak{N}_2(1)     & 
E_1^t=e_1, E_2^t=-e_1+t^{-1}e_2, E_3^t=e_3, E_4^t=-e_3+t^{-1}e_4 \\
\hline

\mathfrak{L}_7 \to  \mathfrak{N}_9(\alpha)     & 
E_1^t=t(\alpha e_2 +e_3), E_2^t=te_1, E_3^t=t^2e_3, E_4^t=t^2e_4 \\
\hline

\mathfrak{L}_7 \to  \mathfrak{L}_1^{\mathbb{C}}    & E_1^t=e_1, E_2^t=e_3, E_3^t=e_4, E_4^t=te_2 \\
\hline
\hline

\mathfrak{L}_8 \to  \mathfrak{N}_2(0)     & 
E_1^t= t(e_1 +e_2-e_3), E_2^t=-t(e_2-e_3), E_3^t=t^2e_3, E_4^t=-t^2e_4 \\
\hline

\mathfrak{L}_8 \to  \mathfrak{N}_3( i)     & 
E_1^t= t(e_1 +e_3), E_2^t= it(e_1-e_3), E_3^t=te_2, E_4^t= t^2e_4 \\
\hline

\mathfrak{L}_8 \to  \mathfrak{L}_7     & 
E_1^t= t^2(e_1 + e_2), E_2^t=t^4(e_2-e_3), E_3^t=t^4(e_3+e_4), E_4^t=t^6e_4 \\
\hline
\hline

\mathfrak{L}_9 \to  \mathfrak{N}_0     & 
E_1^t=te_1, E_2^t=-e_2+\frac{1}{2}e_3-\frac{1}{4}e_4, E_3^t=-e_2-\frac{1}{2}e_3+\frac{1}{4}e_4, E_4^t=-te_4 \\
\hline

\mathfrak{L}_9  \to  \mathfrak{N}_6      & 
E_1^t= te_1, E_2^t= e_2, E_3^t= -t(e_3-e_4), E_4^t= te_3 \\
\hline

\mathfrak{L}_9  \to  \mathfrak{L}_{12}      & 
E_1^t= t^{-1}(e_1+e_2), E_2^t= e_2, E_3^t= t^{-1} e_3, E_4^t= t^{-2}e_4 \\
\hline
\hline

\mathfrak{L}_{10}  \to  \mathfrak{N}_{10}      & 
E_1^t= te_3, E_2^t= te_1, E_3^t= te_2, E_4^t= t^2e_4 \\
\hline

\mathfrak{L}_{10}  \to  \mathfrak{L}_{12}      & 
E_1^t= e_1+te_2, E_2^t= -te_1, E_3^t= t^2 e_3, E_4^t= t^2e_4 \\
\hline
\hline

\mathfrak{L}_{11}  \to  \mathfrak{N}_5      & 
E_1^t= te_1, E_2^t= -te_3, E_3^t= te_2, E_4^t= t^2e_4 \\
\hline

\mathfrak{L}_{11}  \to  \mathfrak{L}_9      & 
E_1^t= 2ie_1-2e_2, E_2^t= -\frac{it}{2}e_2, E_3^t= t(e_3+ie_4), E_4^t= 2ite_4 \\
\hline

\mathfrak{L}_{11}  \to  \mathfrak{L}_{10}      & 
E_1^t= t^{-1}e_1, E_2^t= t^{-2}e_2, E_3^t= t^{-3} e_3, E_4^t= t^{-4}e_4 \\
\hline
\hline

\mathfrak{L}_{12}  \to  \mathfrak{N}_9(-1)      & 
E_1^t= te_2, E_2^t= e_1, E_3^t= e_4, E_4^t= te_3 \\
\hline

\mathfrak{L}_{12}  \to  \mathfrak{L}_1      & 
E_1^t= -t^{-1}e_1, E_2^t= t^{-2}e_2, E_3^t= t^{-3} e_3, E_4^t= t^{-4}e_4 \\
\hline

\end{array}$$
\end{center}
\normalsize

Let now prove primary non-degenerations. First of all, note that $A\not\to B$ if either $A\in\mathfrak{Z}$, $B\in\mathfrak{L}\setminus\mathfrak{N}$ or $A\in\mathfrak{L}$, $B\in\mathfrak{Z}\setminus\mathfrak{N}$. All remaining primary non-degenerations are proved in Table 3. We fix some $4$-dimensional space $V$ with basis $e_1$, $e_2$, $e_3$, $e_4$. Everywhere in Table 3, $\mathcal{R}$ is some set algebra structures on $V$. Whenever a structure named $\mu$ appears, we denote by $c_{i,j}^k$ ($1\le i,j,k\le 4$) the structure constants of $\mu$ in the basis $e_1$, $e_2$, $e_3$, $e_4$ and write $UW$ ($U,W\subset V$) for the product of subspaces of $V$ with respect to the multiplication $\mu$. 
We use the notation $S_i=\langle e_i,\dots,e_4\rangle,\ i=1,\ldots,4$. If in a proof of $A\not\to B$ the standard structure $\lambda$ representing $A$ doesn't lie in $\mathcal{R}$, we define also a map $g\in GL(V)$ such that $g*\lambda\in \mathcal{R}$.

\begin{center} \footnotesize Table 3. {\it Primary non-degenerations of Zinbiel and nilpotent Leibniz  algebras of dimension $4$.}
$$\begin{array}{|l|c|}
\hline
\mbox{non-degenerations}  &  \mbox{reasons}\\
\hline
\hline

\mathfrak{Z}_{1}^{\mathbb{C}}  \ \not\to  \ \mathfrak{B}, 
\mathfrak{B} \in  \left\{ \begin{array}{c} \mathfrak{N}_1^{\mathbb{C}}, \mathfrak{N}_2^{\mathbb{C}}(\beta \neq -2),\\ \mathfrak{N}_3^{\mathbb{C}}, 
\mathfrak{N}_9(2)\end{array}\right\} &
\begin{array}{c}
\mathcal{R}=\left\{\mu\left| \begin{array}{c}
                        S_2S_2+S_1\circ S_3=0,\\
                        S_1\circ S_2 \subset S_4,
						2c_{1,2}^4=c_{2,1}^4 \end{array}\right.\right\}\,\,\,\begin{array}{c}g(e_1)=e_1,g(e_2)=e_2,\\g(e_3)=e_4,g(e_4)=e_3\end{array}
\end{array}\vspace{0.1cm}\\
\hline
\hline

\mathfrak{Z}_{1}   \ \not\to  \ \mathfrak{B}, \mathfrak{B} \in \left\{ 
\begin{array}{c} \mathfrak{N}_1^2, \mathfrak{N}_2(\gamma\not=1,\frac{1}{9}),\\
\mathfrak{N}_3\left(\alpha\neq\frac{i}{2}\right),  \mathfrak{N}_6, \mathfrak{N}_7, \mathfrak{N}_{10} 
 \end{array} \right\} &
\begin{array}{c} 
\mathcal{R}=\left\{\mu\left| \begin{array}{c} 
                        S_1\circ S_4+S_2\circ S_3=0,
						S_2S_2+S_1\circ S_3\subset S_4,\\S_1\circ S_2\subset S_3,
						S_1S_1\subset S_2, 2c_{1,2}^3=c_{2,1}^3, 3c_{1,3}^4=c_{3,1}^4,\\c_{2,2}^4c_{1,1}^2=c_{3,1}^4c_{1,2}^3,
						c_{2,2}^4c_{1,1}^3=3c_{1,2}^3(2c_{1,2}^4-c_{2,1}^4)\end{array}\right.\right\}
						\vspace{0.1cm}\\
\end{array}\\

\hline
\hline

\mathfrak{Z}_{3}   \ \not\to  \ \mathfrak{B}, \mathfrak{B} \in\left\{\begin{array}{c} \mathfrak{N}_2\left(\gamma\neq 0,1\right),\mathfrak{N}_3\left(\alpha \neq\frac{i}{3}\right),\\\mathfrak{N}_6,
\mathfrak{N}_7, \mathfrak{N}_{10}\end{array}\right\} &
\begin{array}{c} 
\mathcal{R}=\left\{\mu\left| \begin{array}{c} 
                        S_1\circ S_4+S_2\circ S_3=0,
						S_1\circ S_2\subset S_4,\\
						S_1S_1\subset S_3, 2c_{1,3}^4=c_{3,1}^4
						\end{array}\right.\right\}\vspace{0.1cm}\\
\end{array}\\

\hline
\hline

\mathfrak{Z}_5  \ \not\to  \  \mathfrak{B}, 
\mathfrak{B} \in  \left\{    \mathfrak{N}_3(\alpha\neq i),  \mathfrak{N}_{10} \right\}  &
Ann_R(\mathfrak{Z}_5)> Ann_R(\mathfrak{B})\\

\hline

\mathfrak{Z}_5  \ \not\to  \ \mathfrak{B}, 
\mathfrak{B} \in  \left\{  \mathfrak{N}_2(1), \mathfrak{N}_9(\alpha\neq-1) \right\}&
(\mathfrak{Z}_5)^{(+2)}<\mathfrak{B}^{(+2)} \\
\hline
\hline


\mathfrak{N}_1^2 \ \not\to  \   \ \mathfrak{B}, 
\mathfrak{B} \in  \left\{ \mathfrak{N}_1^{\mathbb{C}}, \mathfrak{N}_2^{\mathbb{C}}(\beta), 
\mathfrak{N}_9(\alpha\neq1)
\right\}  & Z(\mathfrak{N}_1^2)> Z(\mathfrak{B}) \\

\hline
\hline

\mathfrak{N}_1\not\to\mathfrak{N}_0&   
Ann(\mathfrak{N}_1)>Ann(\mathfrak{N}_0)
 \\
\hline
\hline

\mathfrak{N}_2  (\gamma \neq 0,1) \not\to   \ \mathfrak{B}, 
\mathfrak{B} \in  \left\{ \mathfrak{N}_1^2, \mathfrak{N}_6, \mathfrak{N}_7 \right\}&

\begin{array}{c} 
\mathcal{R}(\gamma)=\left\{\mu\left| \begin{array}{c}
                         S_1S_1\subset S_3,S_1\circ S_3=0,
						S_1 S_2 \subset S_4,\\
						\frac{\gamma}{1-\gamma}(c_{1,1}^3c_{2,2}^4-c_{1,2}^4c_{2,1}^3)^2= (c_{2,1}^3)^2(c_{1,1}^4c_{2,2}^4-c_{1,2}^4c_{2,1}^4)\\
						+c_{2,1}^3(c_{1,2}^4-c_{2,1}^4)(c_{1,1}^3c_{2,2}^4-c_{1,2}^4c_{2,1}^3)
						 \end{array}\right.\right\} \vspace{0.1cm}\\
\end{array}\\

\hline

\mathfrak{N}_2  (\gamma \neq 0,1) \ \not\to  \   
\mathfrak{B},\mathfrak{B}\in\{\mathfrak{N}_0,\mathfrak{N}_3(\beta)\}
&
 Ann(\mathfrak{N}_2(\gamma)) > Ann( \mathfrak{B}) \\

\hline

\mathfrak{N}_2(1)  \ \not\to  \   \mathfrak{B}, 
\mathfrak{B} \in   \left\{ \begin{array}{c} \mathfrak{N}_0, \mathfrak{N}_1^{\mathbb{C}}, \mathfrak{N}_{2}^{\mathbb{C}}(\beta\neq0),\\
\mathfrak{N}_3^{\mathbb{C}},\mathfrak{N}_9(\alpha)   \end{array} \right\} &
\begin{array}{c} Ann_R(\mathfrak{N}_2(1))> Ann_R( \ \mathfrak{B})  \end{array}\\

\hline
\end{array}$$

\end{center}

\begin{center}\footnotesize
$$\begin{array}{|l|c|}
\hline
\mbox{non-degenerations}  &  \mbox{reasons}\\

\hline
\hline

\mathfrak{N}_3(\alpha) \ \not\to  \   \mathfrak{B}, 
\mathfrak{B} \in  \left\{ \mathfrak{N}_2(1), \mathfrak{N}_9 (\alpha)
\right\}   &
\begin{array}{c} \big(\mathfrak{N}_3(\alpha)\big)^2 < \mathfrak{B}^2 \end{array}\\
\hline

\mathfrak{N}_3(\alpha\neq 0) \ \not\to  \  \ \mathfrak{B}, 
\mathfrak{B} \in  \left\{  \mathfrak{N}_3(0), \mathfrak{N}_{10} \right\}&
\begin{array}{c} 
\mathcal{R}(\alpha)=\left\{\mu\left| \begin{array}{c}
                       S_1S_1\subset S_4, 
                        S_4\circ A=0,
                        e_3\in Z(\mu),\\
            \left(\left(\frac{c_{1,2}^4+c_{2,1}^4}{2}\right)^2+\left(\frac{c_{1,2}^4-
            c_{2,1}^4}{2\alpha}\right)^2\right)c_{3,3}^4
            -c_{1,1}^4c_{2,2}^4c_{3,3}^4\\
            -(c_{1,2}^4+c_{2,1}^4)c_{1,3}^4c_{2,3}^4+c_{1,1}^4(c_{2,3}^4)^2+c_{2,2}^4(c_{1,3}^4)^2=0
						 \end{array}\right.\right\} \vspace{0.1cm}\\
\end{array}\\
\hline

\mathfrak{N}_3(0)  \ \not\to  \  \ \mathfrak{B}, \mathfrak{B} \in  \left\{  \mathfrak{N}_{1}^{\mathbb{C}},
\mathfrak{N}_{2}^{\mathbb{C}}(\beta)\right\} 
& Z(\mathfrak{N}_3(0))>Z(\mathfrak{B}) \\

\hline
\hline

\mathfrak{N}_4 \ \not\to  \   \  \mathfrak{N}_{10}   & 
\mathcal{R}=\left\{\mu\left| \begin{array}{c} S_1S_1 \subset S_4,
S_2S_3+S_3S_2 =0,
c_{1,3}^4=c_{3,1}^4\\
	 \end{array}\right.\right\} \\

\hline

\mathfrak{N}_4 \ \not\to  \   \ \mathfrak{B}, 
\mathfrak{B} \in  \left\{\mathfrak{N}_2(1), \mathfrak{N}_9(\alpha)\right\}   &
\begin{array}{c} (\mathfrak{N}_4)^2 < \mathfrak{B}^2  \end{array}\\
\hline
\hline

\mathfrak{N}_5  \ \not\to  \  
\mathfrak{N}_9(\alpha)
&
\begin{array}{c} (\mathfrak{N}_5)^2 < \big(\mathfrak{\mathfrak{N}_9(\alpha)}\big)^2 \end{array}\\
\hline
\hline

\mathfrak{N}_7\not\to\mathfrak{N}_0&   
Ann(\mathfrak{N}_7)>Ann(\mathfrak{N}_0)
 \\
\hline
\hline

\phantom{\!\!\!}\begin{array}{l}
\mathfrak{N}_9 (\alpha\neq -1)  \ \not\to  \   \mathfrak{N}_{1}^{\mathbb{C}}\\ 
\mathfrak{N}_9(\alpha)  \ \not\to  \   \mathfrak{N}_{2}^{\mathbb{C}}(\beta),(1-\alpha)^2\beta+\alpha\not=0\\ 
\mathfrak{N}_9(\alpha\neq 1)  \ \not\to  \   \mathfrak{N}_{3}^{\mathbb{C}}
\end{array}
& 
\begin{array}{c} 
\mathcal{R}(\alpha)=\left\{\mu\left| \begin{array}{c}
                        S_1S_1 \subset S_3, (S_2)^2=0,
                        S_1\circ S_3=0,\\ 
                        S_1\circ S_2 \subset S_4,
						 c_{1,2}^4 = \alpha c_{2,1}^4 \end{array}\right.\right\}\\\vspace{0.1cm}
						 \begin{array}{c}g(e_1)=e_2,g(e_2)=e_1,g(e_3)=e_3,g(e_4)=e_4\end{array}
\end{array}\\
\hline
\hline

\mathfrak{L}_{1}  \ \not\to  \  \mathfrak{N}_{1}^{\mathbb{C}^2}    & 
 (\mathfrak{L}_{1})^{(+2)} < ({\mathfrak{N}_{1}^{\mathbb{C}^2}})^{(+2)}

\\

\hline
\hline

\mathfrak{L}_2 \ \not\to   \ \mathfrak{B}, 
\mathfrak{B} \in   
\left\{ 
\mathfrak{N}_0, \mathfrak{N}_1^{\mathbb{C}},  \mathfrak{N}_2^{\mathbb{C}}(\beta\neq 0), \mathfrak{N}_3^{\mathbb{C}}, \mathfrak{N}_2(1)
\right\} & 
Ann_L(\mathfrak{L}_2 )>Ann_L( \mathfrak{B})  \\

\hline
\hline

\mathfrak{L}_3  \ \not\to   \ \mathfrak{B}, 
\mathfrak{B} \in  \left\{    \mathfrak{N}_1^2,  \mathfrak{N}_2(\gamma\neq 1), \mathfrak{N}_3(\alpha), \mathfrak{N}_7
\right\}  & 
msub_0(\mathfrak{L}_3)>msub_0(\mathfrak{B} ) \\

\hline
\hline

\mathfrak{L}_5  \ \not\to   \ \mathfrak{B}, 
\mathfrak{B} \in  \left\{  \mathfrak{N}_3(\alpha\neq i), 
\mathfrak{N}_{10}, \mathfrak{L}_1\right\}   &
\begin{array}{c}
Ann_L(\mathfrak{L}_{5})>Ann_L( \mathfrak{B}) 
\end{array}\\
\hline

\mathfrak{L}_5  \ \not\to  \ \mathfrak{L}_2  & 
(\mathfrak{L}_5)^2 < (\mathfrak{L}_2)^2 \\
\hline
\hline

\mathfrak{L}_6  \ \not\to   \ \mathfrak{B}, 
\mathfrak{B} \in  \left\{  \mathfrak{N}_1^2, \mathfrak{N}_2(\gamma\neq 1), \mathfrak{N}_6 \right\}& 
\begin{array}{c} 
\mathcal{R}=\left\{\mu\left| \begin{array}{c}
                        S_1S_1 \subset S_3, S_1S_3+S_4 S_1+S_3S_2=0, \\
                        S_1S_2+S_3S_1  \subset S_4, c_{1,1}^3c_{2,2}^4=c_{1,2}^4c_{2,1}^3\\
						 \end{array}\right.\right\} \vspace{0.1cm}\\
\end{array}\\
\hline
\hline

\mathfrak{L}_8  \ \not\to  \  \ \mathfrak{B}, 
\mathfrak{B} \in  \left\{  \mathfrak{N}_2\left(\gamma\neq 0,1\right),\mathfrak{N}_6,\mathfrak{N}_7,\mathfrak{L}_4 \right\} &

\begin{array}{c} 
\mathcal{R}=\left\{\mu\left| 
                        S_1\circ S_2 \subset S_4,
						 \right.\right\} \vspace{0.1cm}\\
\end{array}\\

\hline
\hline

\mathfrak{L}_9  \ \not\to  \  \mathfrak{N}_{10}    & 
msub_0(\mathfrak{L}_9)>msub_0(\mathfrak{N}_{10})\\

\hline
\hline

\mathfrak{L}_{10}  \ \not\to  \  \ \mathfrak{B}, 
\mathfrak{B} \in  \left\{ \mathfrak{N}_0,  \mathfrak{N}_2^{\mathbb{C}}\left( \beta\neq\frac{1}{4} \right),\mathfrak{N}_3^{\mathbb{C}} \right\}    & 
AZ(\mathfrak{L}_{10})>AZ(\mathfrak{B}) \\

\hline
\hline

\mathfrak{L}_{11}  \ \not\to  \  \mathfrak{N}_{3}(\alpha)     & 
AZ(\mathfrak{L}_{11})>AZ\big(\mathfrak{N}_{3}(\alpha)\big)\\
\hline

\mathfrak{L}_{11}  \ \not\to  \ \mathfrak{B}, 
\mathfrak{B} \in  \left\{ \mathfrak{N}_{8},\mathfrak{N}_{9}(\alpha \neq -1) \right\}   & 
(\mathfrak{L}_{11})^{(+2)} <  \mathfrak{B}^{(+2)}\\
\hline

\end{array}$$
\end{center}

\end{Proof}

\section{Irreducible components and rigid algebras}

Using Theorem \ref{theorem}, we describe the irreducible components and the rigid algebras in $\mathfrak{N}$, $\mathfrak{Z}$, and  $\mathfrak{L}$.

\begin{corollary}\label{ir_LZ} The irreducible components  of
$\mathfrak{N}$ are
$$
\begin{aligned}
\mathcal{C}_1&=\overline{\{ O(\mathfrak{N}_2(\gamma)) \}_{\gamma\in\mathbb{C}} }=
\mathfrak{N} \setminus\{  \mathfrak{N}_3(\alpha), \mathfrak{N}_4, \mathfrak{N}_5, \mathfrak{N}_{10} \}_{\alpha\in\mathbb{C}};\\
\mathcal{C}_2&=\overline{ \{O(\mathfrak{N}_3(\alpha)) \}_{\alpha\in\mathbb{C}}}=
 \{ \mathfrak{N}_3(\alpha), \mathfrak{N}_4, \mathfrak{N}_5, \mathfrak{N}_{10},  \mathfrak{N}_1^{\mathbb{C}^2}, \mathfrak{N}_1^{\mathbb{C}}, \mathfrak{N}_2^{\mathbb{C}}(\beta), \mathfrak{N}_3^{\mathbb{C}}, \mathbb{C}^4\}_{\alpha,\beta \in\mathbb{C}}.
 \end{aligned}
 $$ 
In particular, $Rig(\mathfrak{N})= \varnothing.$
\end{corollary}

\begin{Proof} 

The description of irreducible components follows from the proof of Theorem \ref{theorem} and Table 4 given below. The arguments from Table 3 can be applied to prove the required non-degenerations due to Lemma \ref{gmain}.

 \begin{center} \footnotesize Table 4.   {\it Orbit closures for $\mathfrak{N}_2(*)$ and $\mathfrak{N}_3(*)$.}
$$\begin{array}{|l|l|l|}
\hline
\mbox{assertions}  &  \mbox{parametrized bases} & \mbox{parametrized indices} \\
\hline
\hline

\mathfrak{N}_2(*) \to \mathfrak{N}_1 & E_1^t=te_1,E_2^t=te_2,E_3^t=t^2e_4,E_4^t=te_3 & \epsilon(t)=-\frac{1}{t} \\

\hline

\mathfrak{N}_2(*) \to \mathfrak{N}_7 & 
\begin{array}{l}
E_1^t=e_1+(1+t)e_2,E_2^t=-e_1-(1-t)e_2,\\ 
E_3^t=te_3+te_4, E_4^t=-te_3-(t+2t^2)e_4        
\end{array}& \epsilon(t)= 1-t^2\\
\hline

\mathfrak{N}_3(*) \to \mathfrak{N}_4 & E_1^t=t^2e_1+it^2e_3,E_2^t=te_2+t^3e_3,E_3^t=-ie_3,E_4^t=t^2e_4 & \epsilon(t)= \frac{1}{t} \\
\hline

\mathfrak{N}_3(*) \to \mathfrak{N}_5 & E_1^t=e_1,E_2^t=te_2,E_3^t=e_3,E_4^t=e_4 & \epsilon(t)= \frac{1}{t} \\
\hline

\end{array}$$
\end{center}
\end{Proof}

Using Theorem \ref{theorem} and Corollary \ref{ir_LZ}, we obtain the following corollaries.

\begin{corollary} The irreducible components  of 
$\mathfrak{Z}$ are
$$
\begin{aligned}
\mathcal{C}_1&=\overline{\{O(\mathfrak{N}_2(\gamma))\}_{\gamma \in\mathbb{C}}}=
\mathfrak{N} \setminus\{  \mathfrak{N}_3(\alpha), \mathfrak{N}_4, \mathfrak{N}_5, \mathfrak{N}_{10} \}
_{\alpha\in\mathbb{C}};\\
\mathcal{C}_2&=\overline{ \{O(\mathfrak{N}_3(\alpha))\}_{\alpha \in\mathbb{C}}}=
 \{ \mathfrak{N}_3(\alpha), \mathfrak{N}_4, \mathfrak{N}_5,\mathfrak{N}_{10},  \mathfrak{N}_1^{\mathbb{C}^2}, \mathfrak{N}_1^{\mathbb{C}}, \mathfrak{N}_2^{\mathbb{C}}(\beta), \mathfrak{N}_3^{\mathbb{C}},\mathbb{C}^4\}
 _{\alpha,\beta \in\mathbb{C}};\\
\mathcal{C}_3&=\overline{O(\mathfrak{Z}_1)}=
 \{ \mathfrak{Z}_1, \mathfrak{Z}_2, \mathfrak{Z}_{1}^{\mathbb{C}}, 
 \mathfrak{N}_2(1), \mathfrak{N}_3\left(\frac{i}{2}\right), \mathfrak{N}_9(\alpha),
 \mathfrak{N}_1^{\mathbb{C}^2}, \mathfrak{N}_1^{\mathbb{C}}, \mathfrak{N}_2^{\mathbb{C}}(\beta), \mathfrak{N}_3^{\mathbb{C}},\mathbb{C}^4\}_{\alpha,\beta\in\mathbb{C}};\\
\mathcal{C}_4&=\overline{O(\mathfrak{Z}_3)}=
 \{ \mathfrak{Z}_2, \mathfrak{Z}_3, \mathfrak{Z}_{1}^{\mathbb{C}}, 
 \mathfrak{N}_2(0), \mathfrak{N}_2(1), \mathfrak{N}_3\left(\frac{i}{3}\right), \mathfrak{N}_9(\alpha),
 \mathfrak{N}_1^{\mathbb{C}^2}, \mathfrak{N}_1^2, \mathfrak{N}_1^{\mathbb{C}}, \mathfrak{N}_2^{\mathbb{C}}(\beta), \mathfrak{N}_3^{\mathbb{C}},\mathbb{C}^4\}_{\alpha,\beta\in\mathbb{C}};\\
\mathcal{C}_5&=\overline{O(\mathfrak{Z}_5)}=
 \{ \mathfrak{Z}_4, \mathfrak{Z}_5,  
 \mathfrak{N}_3(i), \mathfrak{N}_6, \mathfrak{N}_9(-1),
  \mathfrak{N}_1^{\mathbb{C}^2},  \mathfrak{N}_1^{\mathbb{C}}, \mathfrak{N}_2^{\mathbb{C}}(\beta), \mathfrak{N}_3^{\mathbb{C}},\mathbb{C}^4\}_{\beta\in\mathbb{C}}.
  \end{aligned}
  $$
In particular, $Rig(\mathfrak{Z})= \{ \mathfrak{Z}_1, \mathfrak{Z}_3, \mathfrak{Z}_5 \}$.
\end{corollary}

\begin{corollary}\label{ir_Lei} The irreducible components  of 
$\mathfrak{L}$ are
$$
\begin{aligned}
\mathcal{C}_1&=\overline{\{O(\mathfrak{N}_3(\alpha))\}_{\alpha\in\mathbb{C}}}=
 \{ \mathfrak{N}_3(\alpha), \mathfrak{N}_4, \mathfrak{N}_5,\mathfrak{N}_{10},  \mathfrak{N}_1^{\mathbb{C}^2}, \mathfrak{N}_1^{\mathbb{C}}, \mathfrak{N}_2^{\mathbb{C}}(\beta), \mathfrak{N}_3^{\mathbb{C}},\mathbb{C}^4\}_{\alpha,\beta\in\mathbb{C}};\\
 \mathcal{C}_2&=\overline{ O(\mathfrak{L}_2)}=
 \{ \mathfrak{N}_9(0),\mathfrak{N}_1^{\mathbb{C}^2}, \mathfrak{N}_2^{\mathbb{C}}(0), \mathfrak{L}_1^{\mathbb{C}}, \mathfrak{L}_2, \mathfrak{L}_4, \mathbb{C}^4\};\\
\mathcal{C}_3&=\overline{ O(\mathfrak{L}_5)}=
 \mathfrak{L}  \setminus 
 \{ \mathfrak{N}_3(\alpha), \mathfrak{N}_4, \mathfrak{N}_5,\mathfrak{N}_{10}, \mathfrak{L}_1, \mathfrak{L}_2, \mathfrak{L}_9, \mathfrak{L}_{10},  \mathfrak{L}_{11}, \mathfrak{L}_{12},\mathbb{C}^4\}_{\alpha \in\mathbb{C}\setminus\{i\}};\\
\mathcal{C}_4&=\overline{ O(\mathfrak{L}_{11})}=
 \{ \mathfrak{N}_5, \mathfrak{N}_6, \mathfrak{N}_9(-1), \mathfrak{N}_{10}, \mathfrak{N}_1^{\mathbb{C}^2}, \mathfrak{N}_1^{\mathbb{C}}, \mathfrak{N}_2^{\mathbb{C}}(\beta), \mathfrak{N}_3^{\mathbb{C}}, \mathfrak{L}_1,   \mathfrak{L}_9, \mathfrak{L}_{10}, \mathfrak{L}_{11}, \mathfrak{L}_{12},\mathbb{C}^4\}_{\beta\in\mathbb{C}}.
 \end{aligned}
 $$
In particular, $Rig(\mathfrak{L})= \{ \mathfrak{L}_2, \mathfrak{L}_5, \mathfrak{L}_{11} \}$.
\end{corollary}

\end{document}